\documentclass[twocolumn,journal]{IEEEtran}
\usepackage{algorithm} 
\usepackage{algorithmic} 
\usepackage{authblk}
\usepackage[svgnames]{xcolor}
\usepackage{graphicx} 
\usepackage{subfig} %
\usepackage{caption}
\usepackage{amsfonts, amsmath, amsthm, amssymb} 
\usepackage{multicol}
\usepackage{wrapfig} 
\usepackage{mathptmx}
\usepackage{lettrine}
\usepackage{hyperref} 
\usepackage{cleveref}
\hypersetup{
	colorlinks=true,
	linkcolor=blue,
	filecolor=blue,
	citecolor = blue,      
	urlcolor=cyan,
}
\usepackage[multiple]{footmisc}
\newcommand\tab[1][1cm]{\hspace*{#1}}

\def\keywords#1{\begin{flushleft}{\bf Keywords}\\{#1}\end{flushleft}} 
\date{\today}
	
\begin{document}
\title{ Geodesic Fiber Tracking in White Matter using Activation Function}
\author[1]{Temesgen Bihonegn}
\author[1,3]{Sumit Kaushik} 
\author[1]{Avinash Bansal}
\author[2]{Lubom\'{\i}r Vojt\'\i\v sek}
\author[1]{Jan Slov\'ak}
\affil[1]{Department of Mathematics and Statistics, Masaryk University}
\affil[2]{Brain and Mind Research Programme, Central European Institute of Technology, Masaryk University}
\affil[3]{Faculty of Electrical Engineering, Czech Technical University in Prague, Czech Republic}
\maketitle

\begin{abstract}
\emph{Background and objective:} 
The geodesic ray-tracing method has shown its effectiveness for the
reconstruction of fibers in white matter structure. Based on reasonable metrics
on the spaces of the diffusion tensors, it can provide multiple
solutions and get robust to noise and curvatures of fibers. The choice of
the metric on the spaces of diffusion tensors has a significant 
impact on the outcome of this method. 
Our objective is to suggest metrics and modifications of the algorithms
leading to more satisfactory results in the construction of white matter
tracts as geodesics. \\
\emph{Methods:}
Starting with the DTI modality, we propose to rescale the 
initially chosen metric on the space of 
diffusion tensors to increase the geodetic cost in the isotropic regions.
This change should be conformal in order to preserve the angles between
crossing fibers. We also suggest to enhance the methods to be more robust to
noise and to employ the fourth order tensor data in order to handle the fiber
crossings properly.   
\\  
\emph{Results:} We propose a way to choose the appropriate conformal
class of metrics where the metric gets scaled according to tensor
anisotropy. We use the logistic functions, which are commonly used in
statistics as cumulative distribution functions. To prevent deviation of
geodesics from the actual paths, we propose a hybrid ray-tracing approach. 
Furthermore, we suggest how to employ diagonal projections of 4th order 
tensors to perform
fiber tracking in crossing regions.  
\\
\emph{Conclusions:} The algorithms based on the newly suggested methods were
succesfuly implemented, their performance was tested on both synthetic and
real data, and compared to some of the
previously known approaches. 
\end{abstract} 

	
\section{Introduction}

DTI (Diffusion Tensor Imaging) has become a clinical standard for studying
and diagnosing neuro diseases.  It is the non-invasive approach to
obtain information on the neural architecture.  Fiber tracking
methods broadly comprise of two classes, probabilistic
\cite{TractographyLAZAR2005524} \cite{Tractographyparker}
\cite{probablistic}, and deterministic \cite{Deterministicbasser}
\cite{deterministicConturo10422}
\cite{deterministicLazar2003WhiteMT}
\cite{deterministicdoi:10.1002/1531-8249(199902)45:2<265::AID-ANA21>3.0.CO;2-3}
\cite{deterministicWESTIN200293}.  Probabilistic fiber tracking
traverses all possible trajectories and provides a simulated
distribution of the fiber tracts, which can be used in brain
connectivity studies.  Deterministic tractography methods are
primarily based upon streamline algorithms where the local tract
direction is defined by the principal eigenvector of the diffusion
tensor.  These approaches have been used to construct white matter
anatomical connections in the human brain.  In this work, we are
considering the latter class.

Earlier classical streamline based techniques  \cite{1384} showed
ineffectiveness in the reconstruction of highly curved fibers.  Other
difficulties with these methods appear in the isotropic (slightly
anisotropic) regions where direction information is redundant.  Apart from
that, these methods are also sensitive to noise and fiber crossings.

To overcome the problems mentioned above
\cite{lenglet200410.1007/978-3-540-24673-2_11} \cite{inproceedings}
\cite{Newapproach} \cite{finsler} propose methods based on geodesics in
Riemannian geometric space.  These geodesics follow the shortest path
locally between two points lying on the manifold.  This path is optimal for
the underlying actual fiber tracts.  One class of such methods is based on
Hamilton-Jacobi (HJ) formalism.  These methods are sensitive to local
changes and provide a single solution.

In the works \cite{lenglet200410.1007/978-3-540-24673-2_11}
\cite{Newapproach} authors proposed to use the inverse of diffusion tensor
as the metric tensor of the geometric space.  Fuster and others
\cite{Fuster2015AdjugateDT} introduce modification of inverse metric, called
adjugate tensor, which better explained Brownian motion on Riemannian space
and overcame the issue with inverse diffusion tensor.

In \cite{Sepasian2011MultivaluedGT} \cite{Sepasian2012MultivaluedGR} \cite{Sepasian2008ART}  Sepasian et al comes up with multi-valued ray-tracing method for anisotropic medium. Ray-tracing
methods are based on the assumption that, locally in the medium, 
a wave or particle follows a path corresponding to the least action.
% with varying propagation of velocities. 
Consequently, the directions of the path vary.  
These methods are capable of producing multiple geodesics between point and region in the
medium. 

Local variations of geodesics from underlying fibers are taken under
consideration using Euler-Lagrange equations, but while traversing they
deviate from the actual underlying path.  The conformal rescaling or
adaptive Riemannian metric is chosen for tractography in
\cite{adaptiveremannian} and segmentation in their subsequent work
\cite{HAO2014161}.  Similar to their work, authors in
\cite{novelremannianbyfuster} evaluated adjugate instead of $D^{-1}$ with or
without sharpening.  The choice of adjugate tensor as a metric does not
resolve minimizing the Riemannian cost in all anisotropic or nearly
anisotropic regions.

In this work, the contributions are as follows:  
\begin{itemize}
	\item[1.] Starting with the second order tensor data, 
	we present a method to choose the appropriate conformal
	class of metrics where the metric gets scaled according to
	tensor anisotropy. We use the idea that the rotational
	information is related to the anisotropy of the tensor, 
	and logistic function can be exploited to capture it. 
	In particular, the rotational information is misleading in 
	nearly isotropic regions in the presence of noise. 
	The metric tensor is rescaled, according to this information. We
	compare various scalar anisotropies under the activation function. 
	\item[2.] Ray-tracing method deviates from
	the geodesics path in general. This problem is countered by feeding
	back the principal eigenvector direction of underlying
	interpolated tensor to ODE solver. This also enables
	the hybrid ray-tracing method to perform better in high
	curvature regions. We also enhance the method by local interpolation
	based on the so-called spectral quaternionic distance measures on
	the metric tensors.
	\item[3.] We suggest to employ diagonal projection of 4th order 
	tensors to perform fiber tracking in crossing regions. 
	The diagonal components of the flattened 4th order tensor 
	are second order tensors and lie in Riemannian space. 
	We show that these components have potential to resolve 
	fiber crossings even at small angle intersections. 
%	\item [4.]We need to interpolate; interpolation via spectral metrics works better (because of preserving the
%	anisotropy). We show that spectral quaternion (SQ)/spherical linear version of spectral quaternion (Slerp-SQ) metric based interpolation scheme are better than Euclidean and Log-Euclidean ones.    
	\end{itemize}
	
This paper is organized as follows. In Section II, we review the
geodesic-based fiber tracking approach.  In Section III, parts A through D,
we introduce a modified ray-tracing method, which enables us to find
multiple geodesics by shooting rays from point to region. We describe the
use of activation function, 
%with diagonal component approach for crossing
%fibers by rescaling the classical metric, 
which we call $\beta$-scaled
metric tensor. Next, we suggest how to employ the so called diagonal
components of the fourth order tensors, \cite{SumitKaushik}, to resolve the
fiber crossing even at small angles. Finally, in part G we 
also comment on various choices of metrics suitable for
local interpolation of tensor data.  Section IV shows the results of our
tracking approach on synthetic and real brain diffusion data.
	
\section{Background} 
	In geodesic ray-tracing, a small deviation of the geodesics
	from the direction of diffusion is preferred. It makes the
	geodesics robust to noise, but if this deviation is big, it needs
	a sharpening of diffusion tensor \cite{Sepasian2011MultivaluedGT} \cite{Descoteaux}. It can be done by powering the tensor. However, it causes artificial increase or decrease in volume, which is not required as the diffusion process is physical, and diffusion quality must be preserved. This is partially done by the normalization of the tensor. The
	sharpening strategy seems to result in better tractography. For
	more details, we refer the readers to \cite{Descoteaux}. 
	
	 The works \cite{jbabdi2008accurate} \cite{ThomasFletcher:2013:GRT:2559054.2559117} explain the choice of inverse diffusion tensor as a metric in the context of DTI. It does not work for all cases. Another approach for modification of metric has been indicated by Hao \cite{HAO2014161}, which has a similar effect as the adjugate metric proposed by Fuster \cite{Fuster2015AdjugateDT}. The two approaches are build upon the conformal rescaling of the tensor. They use adjugate tensor with sharpening to track high curvature fibers. The main idea  to use inverse diffusion tensor as the metric tensor is to ensure that path is shorter if diffusion is stronger along the high anisotropic direction. This provides the minimization of the path, in essence, which can be treated as a geodesic. The ray-tracing method works under the assumption to consider a bundle of the rays together and provides a multi-valued solution. In this work, we consider a cone formed on the base of the ellipsoid. 

	Sepasian et al introduced a modified ray-tracing by
	adjusting the direction of geodesics based on computing
	Ricci curvature tensor from the metric tensor and its derivatives \cite{Sepasian2016ModifiedGR}. It
	provides a measure of the degree of deviation determined by the Riemannian metric tensor from Euclidean space. DTI model fails in the regions where fibers are	merging, intersecting, and kissing. 
	The second order tensor in DTI lies in the Riemannian space, which is well studied
	in \cite{ThomasFletcher:2013:GRT:2559054.2559117} \cite{Pennec2006} \cite{lengletremannian} \cite{Fletcher:2003:SSV:1965841.1965853} \cite{Krajsek2016}. 

	The geodesic methods employing
	the Hamilton-Jacobi equation (HJ) fail in highly curved regions
	comparing to fast marching techniques. In  \cite{controlandfast}, geodesics
	are considered as a function of position and direction. In
	isotropic regions, the rays may deviate from the actual
	path \cite{Descoteaux} \cite{Lazar2003WhiteMT} \cite{adaptiveenhancement}. Sharpening is helpful in these cases, as mentioned above.

	During traversal, geodesic rays tend to deviate so that there is
	non-uniformity in their distribution across the regions, which
	causes less dense fibers. Ray density can be altered by
	changing the mesh size of interpolation. Anisotropic diffusion
	in Euclidean space is similar to the Brownian motion of water
	molecules in an isotropic medium in Riemannian Space \cite{Fuster2015AdjugateDT}.
	
\section{Ray-tracing via Activation function}
\subsection{Initial Shooting Direction}

\begin{figure}[htbp!]
	\centering
	\subfloat[elliptic cone]{\includegraphics[width=0.2\textwidth]{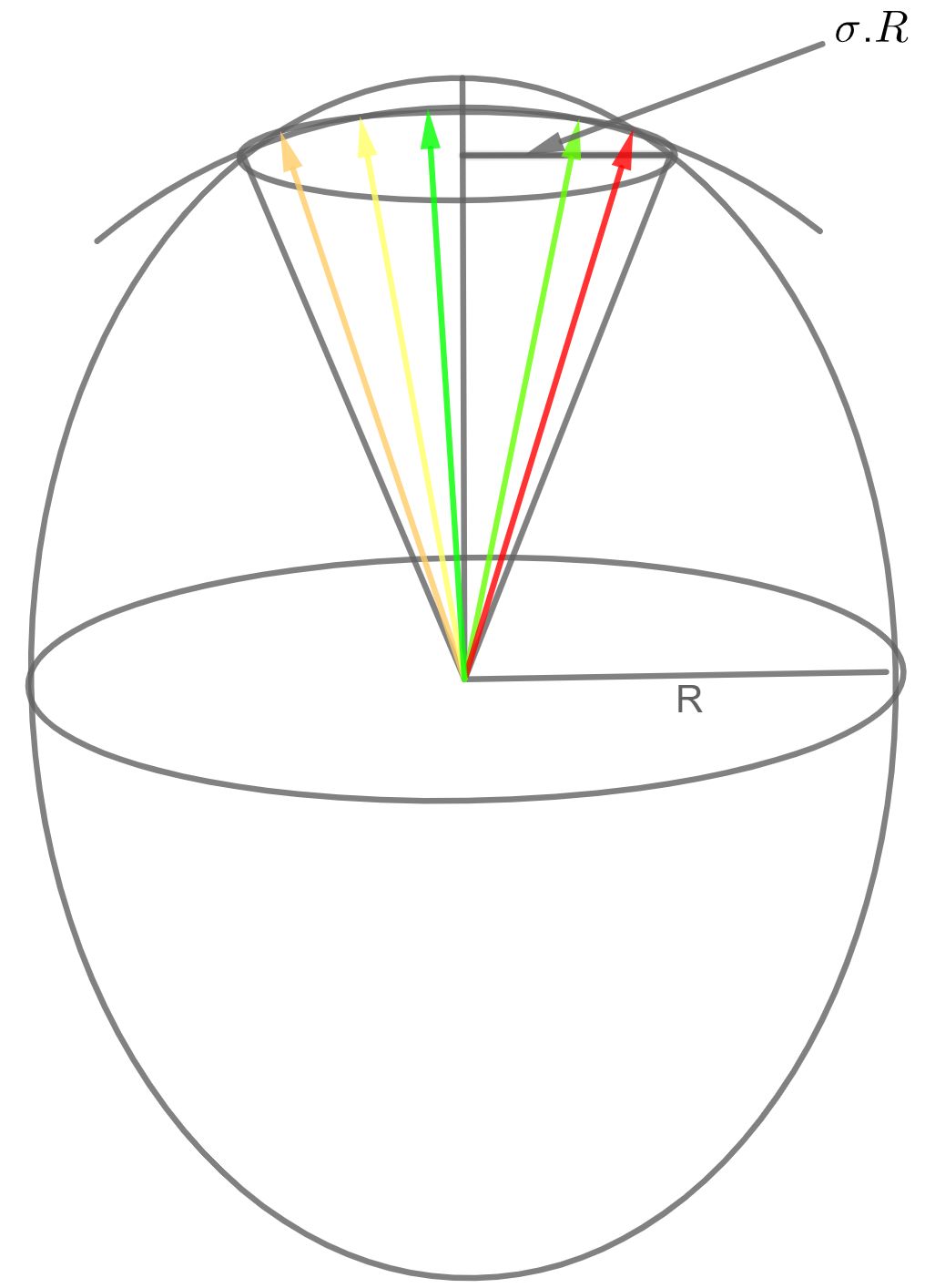}\label{fig:elliptic}}
	\subfloat[cone shooting]{\includegraphics[width=0.3\textwidth]{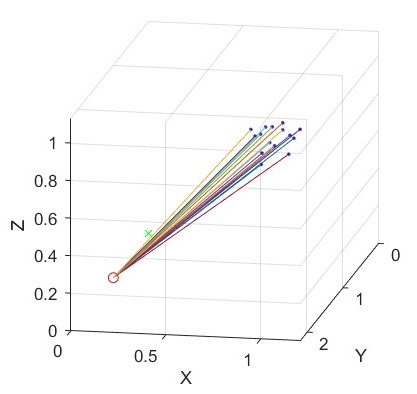}\label{fig:cone}}
	\caption{Initial shooting direction}
	\label{fig:shootingdirection}
	\end{figure}

Ray-tracing method is used to find the trajectory of the particle moving in
medium.  To compute the geodesics using this technique we need the
initial position and direction.  In DTI tractography, we can
restrict the initial shooting directions.  Let $R$ be the radius of
the base and $\sigma\in(0,1)$, which adjusts the base of the cone. 
The different values of $\sigma$ provide different bases
($\sigma.R$) of a cone.  Directions are uniformly distributed over a
spherical section of the cone, as shown in Fig
~\ref{fig:shootingdirection}.  This is done to restrict the shooting
direction and to ensure the ray bundle remains densely packed.

	%\vspace*{-0.3cm}	
	The values of $\sigma$ approaching to 1 causes bigger perturbation. The more realistic way is to ray trace from point to region or region to region because it is not possible to know in advance if the initial and final points are connected \cite{Sepasian2011MultivaluedGT}.

\subsection{Activation Function}
The logistic or activation functions are known for their common use in deep
learning methods and statistics as cumulative distribution functions.  One
of their special cases is sigmoid functions, which are differentiable over
real domain values and have positive derivatives at each point.  To account
for the tensors with negligible difference between the maximal
($\lambda_{max}$) and minimal ($\lambda_{min}$) eigenvalues, a smooth
transition function is applied.\nocite{Collard:2014:APD:2582933.2582989} In
our case, we tested the following three functions with very similar
performance:
%\vspace{0.01cm}
\begin{subequations}\label{eq:group}
	\begin{align}
	 S_{1}(x)=\tanh(x)  \tab[1cm] \label{eq:A}\\
	 S_{2}(x)=\frac{1}{1+\exp(-\frac{1}{2}x)}   \tab[0cm] \label{eq:B}\\
	 S_{3}(x)=\dfrac{x}{\sqrt{1+x^2}}  \tab[0.9cm]\label{eq:C}
	 \end{align}
\end{subequations}
In our experiments discussed below, the function $S_1$ was used.

Hilbert Anisotropy \cite{koufanyhilbert} is given by:
%	\vspace{-1cm}
\begin{equation}\label{eq_2}
HA=\log(\dfrac{\lambda_{max}}{\lambda_{min}}),
\end{equation}
where $HA\geq 0,~HA=0$ for fully isotropic tensor.
 
HA is a scalar measure of anisotropy and is scale-invariant (depends on the shape not the size of the tensor). It is also invariant to rotation and it is a dimensionless number reflecting microscopic diffusion at the level of tissues \cite{Collard:2014:APD:2582933.2582989}. To choose an appropriate metric, we scale the Riemannian metric by an activation function, which is adapted according to the inherent anisotropy.\\
Let $\beta:=S_{i}(x)$  for $ x=HA, i=1,2,3$, then the $\beta$-scaled metric is given by 
%\vspace{0.1cm}
\begin{equation}\label{eq_3}
g_{\beta}=\beta^{-p}D^{-n},
\end{equation}
where $D$ is the second  order diffusion tensor $p\ge 1$, $n\ge 1$. In our experiments, we
used $n=2$, and mostly $p=2$. Recall, $D$
belongs to the space $S^+(3)$ of positive definite $3\times 3$ matrices (SPD). In
particular, $D$ and all its powers are metric tensors.

%When $p<1$ it decreases the Riemannian cost in the isotropic region and
%increase the cost in anisotropic region.  So we choose $p\geq1$ to minimize
%the cost in anisotropic region.

Beside the Hilbert anisotropy, various scalar measures exist, which can serve as description for the degree of anisotropy of diffusion tensor. These measures
can be composed with the above functions. They include: mean diffusivity (MD), fractional anisotropy (FA), relative anisotropy (RA) and geometric ones: geodesics anisotropy (GA), Hilbert anisotropy (HA) \cite{Collard:2014:APD:2582933.2582989}. 

Spectral metrics allow for proper scaling of the rotational contribution according to the anisotropy. This is achieved by using the combination of the activation function with anisotropy scalar measure. In Fig~\ref{fig:effect of interpolation}, we compare the Riemannian cost while considering these anisotropy measures under the activation function. From left to right the tensor exhibit high-low-high spectrum of anisotropies. The interpolation of tensors in between two extreme anisotropic tensors is obtained using Log-Euclidean metric shown in equation \eqref{LogE formula}. Under this metric we can observe that the anisotropy is not preserved. There is variation in eigenvalues and rotational component of interpolated tensors.

The minimal Riemannian cost in anisotropic direction is given by $\beta\lambda_{max}^{-2}$. 
 It is observed that Riemannian cost increases as tensors achieve high isotropy in the middle of the spectrum and after which a smooth descend is noticed for HA case. The other scalar measures do not give linear interpolation as shown in
\cite{Collard:2014:APD:2582933.2582989}, see Fig 3 there. HA is the only one of all above mentioned scalar anisotropy measure keeping affine combinations invariant.

\begin{center}\vspace{.1cm}
	\includegraphics[width=.9\linewidth]{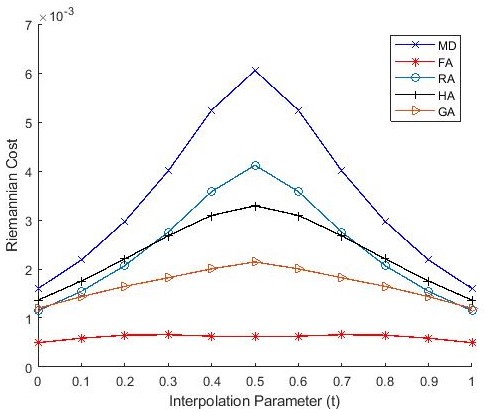}
	\includegraphics[width=.85\linewidth]{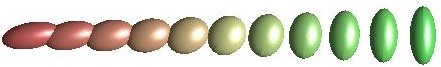}
	\captionof{figure}{\color{Black} The effect of interpolation between the tensors and the Riemannian cost from anisotropy to isotropic region. HA gives low Riemannian cost for anisotropic and  high Riemannian cost for isotropic region.}
	\label{fig:effect of interpolation}
\end{center}\vspace{.01cm}

\subsection{Governing Equations}
The trajectory of a fiber pathway is computed iteratively from the hybrid approach, position from ODE
solver, and direction equal to the principal eigenvector direction. The geodesic method in Riemannian manifold which is used to compute the trajectory of the fibers from ODEs is shown below. 

Let $x(\tau)$ be a smooth and differentiable parametrized curve in the Riemannian manifold, $\tau=[0,T]$. The Riemannian length is given as follows:\\
  % \end{flushleft}
	\vspace{0.1cm}
\begin{equation}\label{eq_4}
L(x,\dot{x})=\int_{0}^{T}(g_{\alpha\beta}\dot{x}^{\alpha}\dot{x}^{\beta})^{1/2}d\tau
\end{equation}
%\begin{flushleft}
The geodesic is the curve that minimizes the length~\eqref{eq_4}. The technique of the Euler-Lagrange equations for solving variational problems is
explained, e.g., in \cite{jost2005riemannian}.

Let $\dot{x}^{\gamma}$ and $\ddot{x}^{\gamma}$ be the first and second derivative with respect to $\tau$, respectively, of the geodesic for dimension $\gamma={1,2,3}$. The geodesics are given by the following system of equations \vspace{0.1cm}
	
\begin{equation}\label{eq_5}
	\ddot{x}^{\gamma}+\sum_{\alpha=1}^{3}\sum_{\beta=1}^{3}\Gamma_{\alpha\beta}^{\gamma}\dot{x}^{\alpha}\dot{x}^{\beta}=0,	
\end{equation}
where $\Gamma_{\alpha\beta}^{\gamma}$ are the so called Christoffel symbol, given by
\begin{equation}\label{eq_6}
\Gamma_{\alpha\beta}^{\gamma}=\dfrac{1}{2}\sum_{\sigma=1}^{3}g^{\gamma\sigma}\left( \dfrac{\partial g_{\beta\sigma}}{\partial x^{\alpha}}+\dfrac{\partial g_{\alpha\sigma}}{\partial x^{\beta}}-\dfrac{\partial g_{\beta\alpha}}{\partial x^{\sigma}}\right) \
\end{equation}
and $g_{\beta\sigma}$ denotes the matrix component of the inverse diffusion tensor, and $g^{\gamma\sigma}$ represents an element of the original diffusion tensor.
We compute the solution of equation~\eqref{eq_5} for the given initial position and multiple initial directions using the standard ODE solvers, such as fourth-order Runge-Kutta method. This gives us a set of geodesics connecting the given initial point, which we integrate until they hit the boundary. Depending on the equations~\eqref{eq_5} and ~\eqref{eq_6}, we need nine symbols per dimension, for a total of 27 symbols. However, dealing with torsion free
connections allows to exploit additional symmetries. The initial position is user-specified and directions are computed by forming a cone with a base of the radius. 

The choice  of power of $\beta$ is done experimentally. In our work, we have compared the results on synthetic data with $D^{-1}$, adjugate and
$\beta$-scaled diffusion tensor. The experiment shows our approach works irrespective of configuration in terms of curvature and (an)isotropy of neighboring tensors.

Based on the observation that the ODE
solver's output direction deviates from the actual fiber path, we used the principal eigenvector of the
underlying interpolated tensor as input for the ODE solver. While
picking up the principal vector direction, there are always
possibilities of choosing two directions. At each iteration, we
need to keep track of following the direction consistent with
traversing fiber. This hybrid approach resulted
in the traversal of geodesics in the high curvature cases and
is robust to noise as well.\\
\vspace{0.1cm} 
\begin{algorithm}  	\label{Hybrid raytracing}
	\floatname{algorithm}{Algorithm 1}
	\renewcommand{\thealgorithm}{}
	\caption{Hybrid Ray-Tracing}

	 \hspace*{\algorithmicindent} \textbf{Input}: Initial position ($x$) and direction ($\dot{x}$) \\
	\hspace*{\algorithmicindent} \textbf{Output}: Local geodesics 
	\begin{algorithmic}[1]
		\STATE Define the mesh size locally over the physical grid with size $m=0.1$  
		\STATE Compute interpolated inverse diffusion tensor locally.
		\STATE Find local geodesics
		\begin{itemize}
			\item[a] Give the position and direction to ODE solver.
			\item[b] Compute Christoffel symbols using Algorithm 2.
		\end{itemize}
		\STATE Take the new position and replace directions with principal eigen vector of the underlying interpolated tensor.
		\STATE Repeat step $2-4$ until the geodesics leave the grid.
	\end{algorithmic}
\end{algorithm} 

\begin{algorithm}
	\floatname{algorithm}{Algorithm 2}
	\renewcommand{\thealgorithm}{}
	\caption{Compute Christoffel symbol $\Gamma^i_{jk}$}
	\label{Christoffel symbol}
	 \hspace*{\algorithmicindent} \textbf{Input}: the diffusion tensor $D$
and the indices $i,j,k$. \\
	\hspace*{\algorithmicindent} \textbf{Output}: Christoffel symbols
	\begin{algorithmic}[1]
		\STATE Adapt $D$ to $g_\beta$ according to ~\eqref{eq_3}
		\STATE Set $\Gamma^i_{jk}=0$ 
		\STATE Loop through all three dimensions, i.e $m=1,2,3$ 
		\STATE   $\Gamma^i_{jk}=\Gamma^{i}_{jk}+\dfrac{1}{2}
\bigl(g_\beta{}^{-1}(m,k).\bigl(D2(j,m,i)+ D2(m,i,j)- \hspace*{1in} D2(i,j,m)\bigr)\bigr)$ 
		\STATE End of loop \\
		Here, D2 is the second order difference of the metric tensor
		$g_\beta$ in \eqref{eq_6}.
	\end{algorithmic}
\end{algorithm} 

\subsection{Rescaling of Metric Tensor}
Illustration of Fig \ref{Anisotropic and isotropic region illustration} comes from \cite{Fuster2015AdjugateDT} that advocate the use of adjugate diffusion tensor instead of the inverse of diffusion tensor as a metric tensor. The intuitive idea is to minimize the Riemannian cost along the trajectory.

\begin{center}\vspace{-0.5cm}
	\includegraphics[width=.8\linewidth]{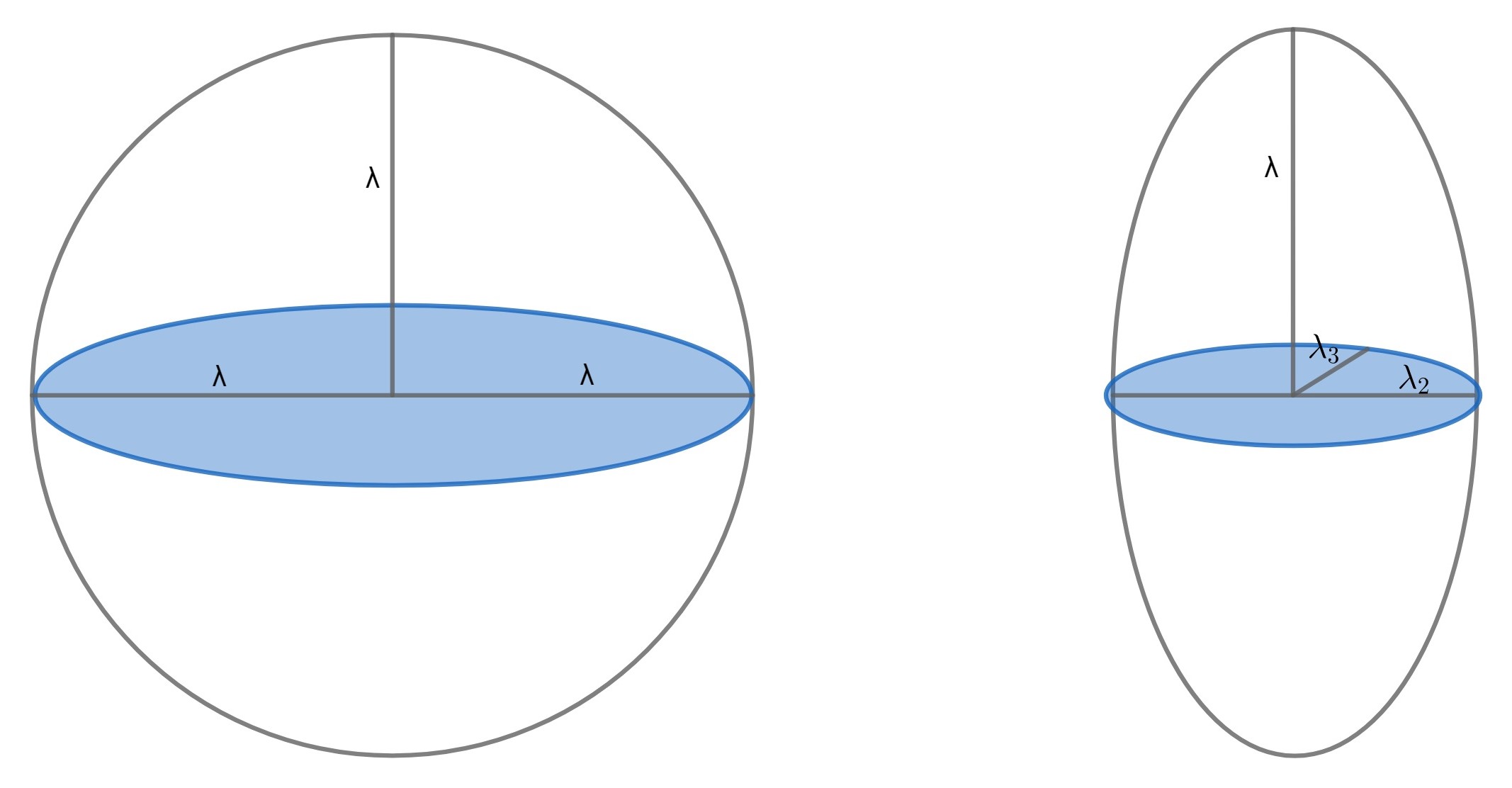}
	\captionof{figure}{\color{Black} An isotropic region (left) and anisotropic region (right)}
	\label{Anisotropic and isotropic region illustration}
\end{center}\vspace{0.1cm} 
Consider two tensors whose principal eigenvalues are equal. In $D^{-1}$ case, Riemannian cost \eqref{eq_4} of (traveling along) an infinitesimal vertical line element scales by $1\setminus\lambda$. For adjugate case i.e., $dD^{-1}$ , where $d=\det(D)$ the Riemannian  cost for isotropic tensor is proportional to  $\lambda^{2}$  (size of the shaded circle) and for anisotropic tensor it is $\lambda_{2}\lambda_{3}$ (proportional to the size of the shaded region).
This method does not work if $\lambda>\lambda_{2},\lambda_{3}$. When the area of the orthogonal cross section in the isotropic
case becomes equal (i.e., same $\lambda_{2},\lambda_{3}$) but their principal eigenvalues are different, adjugate tends to give the same Riemannian cost whereas our approach scales the metric appropriately according to the scalar anisotropy. The scaling coefficient takes zero value for isotropic and higher values for anisotropic cases.

\begin{figure*}[htbp]\label{fig:group}
	\centering
	\subfloat[case1]{\includegraphics[width=0.15\textwidth]{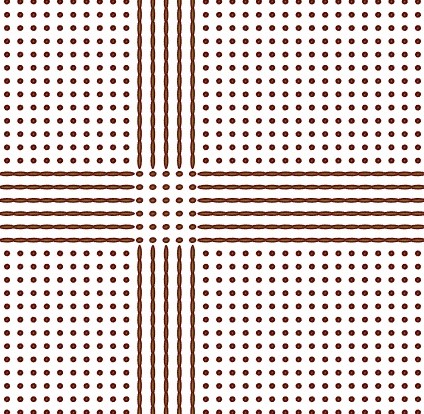}}
	\label{fig:case1}
	\subfloat[$D^{-1}$]{\includegraphics[width=0.18\textwidth]{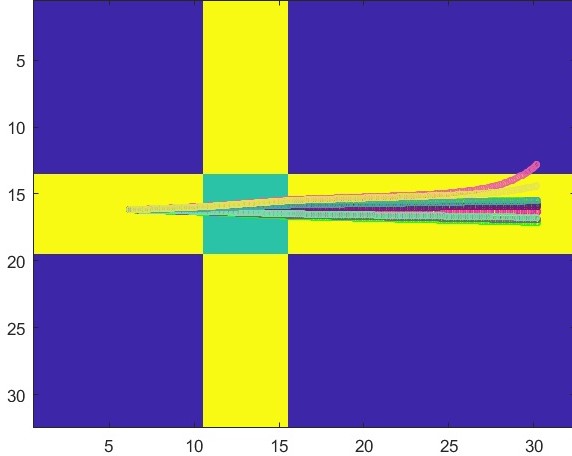}}
	\subfloat[Adjugate]{\includegraphics[width=0.18\textwidth]{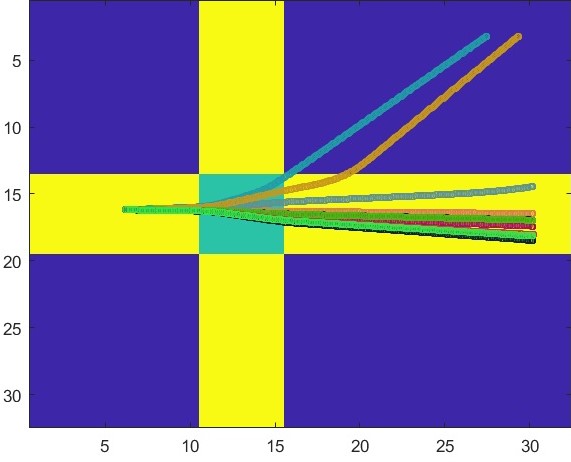}}
	\subfloat[$\beta$-scaled, $p=2$]{\includegraphics[width=0.18\textwidth]{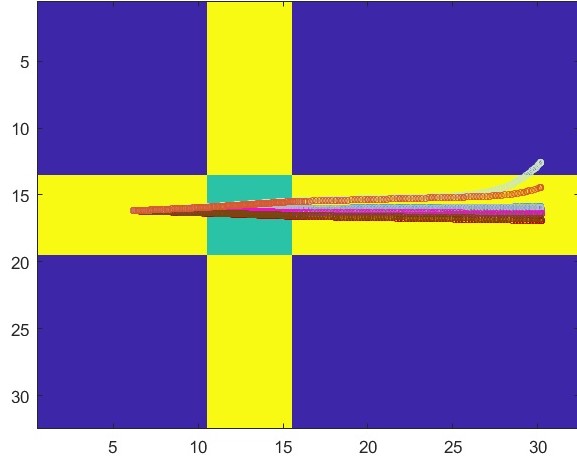}}\\
	\subfloat[case2 ]{\includegraphics[width=0.15\textwidth]{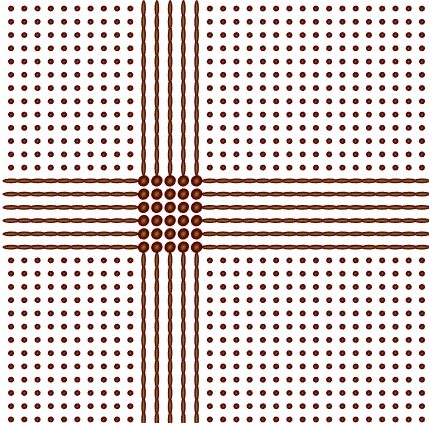}}
	\subfloat[$D^{-1}$]{\includegraphics[width=0.18\textwidth]{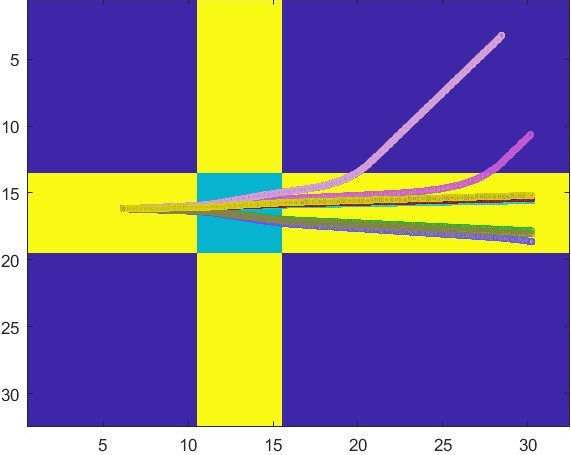}}
	\subfloat[Adjugate]{\includegraphics[width=0.18\textwidth]{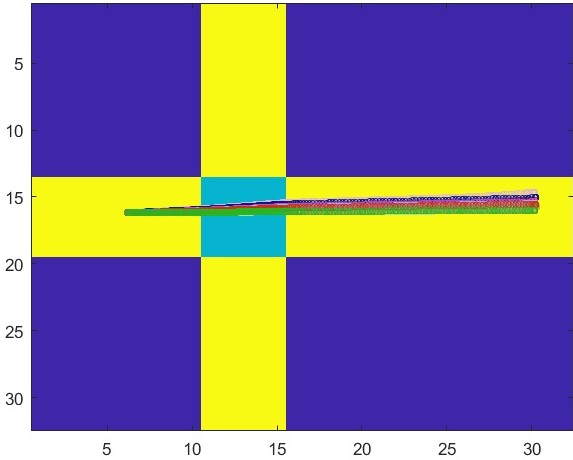}}
	\subfloat[$\beta$-scaled, $p=2$]{\includegraphics[width=0.18\textwidth]{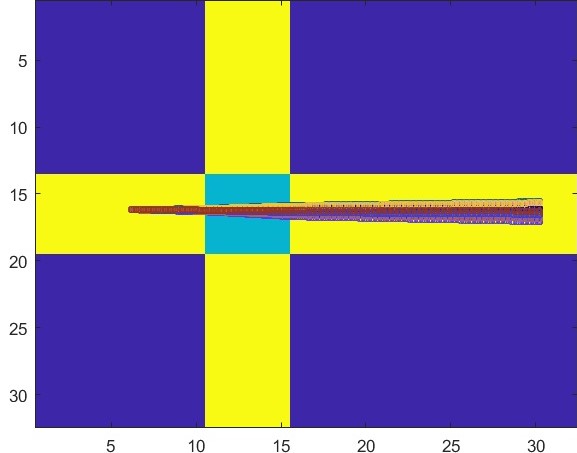}}
	\caption{Comparison of the three metric tensors with the two extreme cases}
	\label{fig:two extreme cases}
\end{figure*}

HA is zero irrespective of the size of the isotropic tensor. This leads to same evaluation of Riemann costs for any isotropic tensor. Such scaling suggests that the diffusion of water molecules is uniform in all directions and hence the Riemann cost as well.
In Fig \ref{fig:two extreme cases} two cases are depicted. Fig \ref{fig:two extreme cases}(a) shows the case 1, where isotropic tensors in the intersection region are chosen  with the smaller eigenvalues. Fig \ref{fig:two extreme cases} (e) shows the case 2 with larger eigenvalues. In both of the cases, the metrics $D^{-1}$ and adjugate induces different Riemann cost. However, in both cases, the
$\beta$-scaled metric lowers the cost of traversing irrespective of the eigenvalues of isotropic tensors. 

\subsection{Decomposition of 4th Order Tensor}
Tuch \cite{tuch2002high} introduced the idea to use mono-exponential model for diffusion of water molecules in the tissues using multiple gradient directions:
\begin{equation}\label{eq_7}
S =  S_{0}\operatorname{exp}(-bD(g))
\end{equation}
For anisotropic diffusion this equation \eqref{eq_7} is linear in the log domain, thus,
\begin{align*}
\log(S)=\log(S_{0})-bD(g)
\end{align*}
where, 
\begin{align*}
D(g) &=  \sum_{i_{1}=1}^{3}\sum_{i_{2}=1}^{3}\sum_{i_{3}=1}^{3}\cdots
\sum_{i_{n}=1}^{3}D_{i_{1}j_{2}i_{3}\cdots i_{n}}g_{i_1}g_{i_2}g_{i_3}\cdots
g_{i_n}
\end{align*}
Here, ${D}_{i_{1}\cdots i_{n}}$ are the coefficients of $n$-th order tensor, while $g_i$ are components of the unit gradient vector $g$, $b$ is the diffusion weighting
factor, and  $S$ and $S_0$ are drop in the signal in presence and absence of diffusion gradients respectively.
Earlier methods based on the least square estimation do not ensure positive diffusion profile. The methods proposed in \cite{barmpoutis2009regularized} \cite{barmpoutis2010unified} ensures positive semi-definiteness of the tensors. 
 We apply flattening of 4th order tensor, which gives $9\times9$ matrix, and eigen-tensors have the potential to reveal actual fiber directions \cite{jayachandra2008fiber}. The diagonal components approach \cite{SumitKaushik} retains geometrical information of the full tensor. The diagonal component of this matrix is symmetric positive definite tensor \cite{SumitKaushik}. In general, $n^{th}$ order tensor
$T^{(n)}$
can be expressed as a matrix of $(n-2)$ order tensors: 
\begin{equation}
T^{(n)} = \begin{pmatrix}
T_{xx}^{(n-2)} & T_{xy}^{(n-2)} & T_{xz}^{(n-2)}\\
T_{yx}^{(n-2)} & T_{yy}^{(n-2)} & T_{yz}^{(n-2)}\\
T_{zx}^{(n-2)} & T_{zy}^{(n-2)} & T_{zz}^{(n-2)}
\end{pmatrix}
\end{equation}
For instance, the diagonal block element $T_{xx}^{(2)}$ in the fourth order
diffusion tensor is given by 
\begin{equation} \label{eq:9}
%{T_{xx}^{(2)}} =
\begin{pmatrix} T_{xx(xx)} &{} T_{xx(xy)} &{} T_{xx(xz)}\\ T_{xx(xy)} &{} T_{xx(yy)} &{}
T_{xx(yz)}\\ T_{xx(xz)} &{} T_{xx(yz)} &{} T_{xx(zz)} \end{pmatrix}\\
=\begin{pmatrix} D_{xxxx} &{} D_{xxxy} &{} D_{xxxz}\\ D_{xxxy} &{} D_{xxyy} &{} D_{xxyz}\\ D_{xxxz} &{} D_{xxyz} &{} D_{xxzz} \end{pmatrix}
\end{equation}
Another observation is that the flattening of 4th
order tensor using diagonal components (DC) can potentially reveal the
actual underlying fiber directions. This observation could be quite useful 
in fiber tracking. 

\subsection{Resolution of Fiber Crossings at Fine Angles}

We have shown experimentally that these diagonal components produce small
orientation errors in comparison to the Cartesian tensor fiber orientation
distribution (CT-ODF) method \cite{weldeselassie2012symmetric} and
\cite{weldeselassie2010symmetric}.  In \cite{weldeselassie2012symmetric},
Figure 3, shows a comparison of orientation errors computed from CT-ODF
method versus the other methods: QBI, DOT MOVMF and MOW.  

The other
observation is about fuzziness in finding maxima using the CT-ODF method. 
These maxima provide the direction of underlying fibers.  The maximal of ODF
does not necessarily align with the actual underlying fiber direction.  The
CT-ODF method does the misalignment correction.  The correction involves the
computation of the maxima.  

\begin{figure}[htbp]
	\centering
	\subfloat[angle between fibers $60^{\circ}$]{\includegraphics[width=0.18\textwidth]{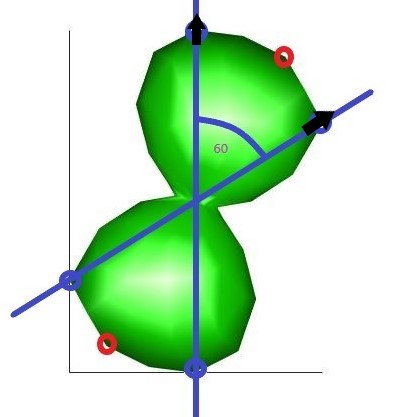}} 
	\subfloat[angle between fibers $75^{\circ}$]{\includegraphics[width=0.18\textwidth]{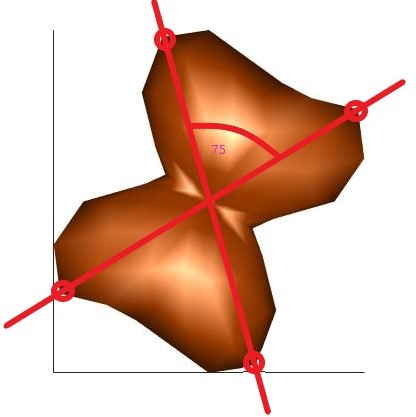}}\
	\caption{4th order ODF with angle differences between the two fibers}
	\label{fig:4th order Using DC}
\end{figure}

The CT-ODF method shows ambiguity in finding
maxima, as shown in Fig \ref{fig:4th order Using DC}.  In Fig~\ref{fig:4th
order Using DC}(a), when the angle between fibers is less than $70^{\circ}$,
the points labeled by arrows are supposed to be the better choices for
maxima than the middle one (point labelled by red colour).  This effect
disappears when the angle difference falls in $70\leq\theta\leq90$ range
(see Fig~\ref{fig:4th order Using DC}(b)).  

\begin{figure}[htbp]
	\centering
	\subfloat{\includegraphics[width=0.5\textwidth]{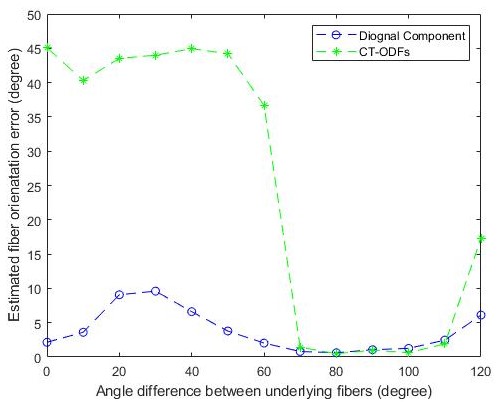}} 
	\caption{Comparison of DC vs. CT-FOD}
	\label{fig:DCvs CTODF}
\end{figure}

Our projection to second order
tensors mentioned in last section is devoid of this ambiguity.  In Figure
\ref{fig:DCvs CTODF}, the known angle differences between the two fibers are
shown on x-axis.  These fibers are used to generate the ODFs.  These angles
are then retrieved using CT-ODF and D-components.  The resulting orientation
errors are shown on y-axis.  For the angle range $70\leq\theta\leq 110$,
performance of both the methods is comparable whereas for smaller angle
differences there is a significant drop of the estimated error 
visible for D-components.  
%In this range, both methods are comparable.  

In the next subsection, we propose to use CT-ODF for
reorientation and diagonal components for tracking fibers, particularly in
crossing regions. This method is extendable to higher-order tensors; for
instance, 6th order tensor has nine diagonal components which could resolve
up to 
nine directions. However, practically more than 3 or 4 fibers per crossing
seldom arise.

\smallskip
\noindent\textbf{Resolving crossing fibers.}
When dealing with the crossing regions, we enhance Algorithm 1 by working in
two layers corresponding to two projections of the fourth order tensor 
to its diagonal components, see Algorithm 3. 

\begin{algorithm}
	\floatname{algorithm}{Algorithm 3}
	\renewcommand{\thealgorithm}{}
	\caption{Reconstruction of fibers}
	\label{Crossing fibers}
	\hspace*{\algorithmicindent} \textbf{Input}: Reoriented 4th order tensor field\\
\hspace*{\algorithmicindent} \textbf{Output}: Fiber reconstruction
	\begin{algorithmic}[1]
		\STATE Flatten the 4th order tensor
		\STATE Extract two layers corresponding to the diagonal components in 2D using
the matrix representation \eqref{eq:9}
		\STATE Shoot the rays from initial point/region using Algorithm 1 in layer 1 and layer 2. 
	\end{algorithmic}
\end{algorithm}

\subsection{Local Interpolation Effect}\label{local Interpolation}
The interpolation step in the
algorithm affects the flow of geodesics while fiber tracking. Aside 
the elementary Euclidean interpolation of the tensor data, there are smarter
choices available, including the 
Log-Euclidean (LogE), Spectral Quaternion (SQ), and spherical version of
spectral quaternionic interpolation (SlerpSQ) \cite{SumitKaushik}, detailed
explanation can be also found in \cite{Collard:2014:APD:2582933.2582989}.
 
\noindent\textbf{Log-Euclidean Interpolation}: In this geometry, the distance between two tensors $T_{1},T_{2}\in S^{+}(3)$ is given by 
\begin{align}
d_{LogE}(\mathbf{T}_1,\mathbf {T}_2) = ||\mathbf\log( {T}_1)-\mathbf \log({T}_2)||
\end{align}
It is based on the fact that the symmetric $3\times3$ matrices are 
diffeomorphic to $S^{+}(3)$ via the exponential mapping.
The interpolation curve between two tensors is the geodesics curve $\gamma_{LogE}:[0,1]\longrightarrow S^{+}(3)$, where the space $S^{+}(3)$ is a convex subset of the Euclidean space $R^{3\times3}$ of $3\times3$ matrices and it is given
for all $0\le t\le 1$ by:
%\begin{align}
\begin{align}
	\gamma_{LogE}(t)=\operatorname{exp}(t\log(T_{1})+(1-t)\log(T_{2})),%\, for \hspace{0.2cm} 0 \leq t\leq 1
	\label{LogE formula}
\end{align}

\noindent\textbf{Spectral Quaternion Interpolation}: The basic idea of spectral metric is to treat eigenvalues and eigenvectors of a
SPD matrix separately. The eigenvalue decomposition of the SPD matrix in spectral geometry is $T=R\varLambda R^{T}$ into a rotation matrix ${R} \in \mathbb SO(3)$ and a diagonal matrix $\varLambda$ containing the
eigenvalues, which provides a natural way of splitting the tensor. Thus using the spectral decomposition of a positive definite matrix, the interpolation curve
is given by the equations
\begin{align}
{S}(t)&= {R}(t) \varLambda (t) {R}(t)^T,\\
{R}(t)&=  {R}_1 \exp (t \log ( {R}_1^T  {R}_2)), \\
\varLambda(t)&=  {R}_1 \exp (t \log ( {R}_1^T {R}_2))
\end{align}
The geometric interpretation of the interpolation curve is a geodesic in the product space of the Lie group defined as 
  $G= {\mathrm {SO}}(3) \times \mathrm {D}^+(3)$, where  $\mathrm {D}^+(3)$ is the group of diagonal matrices with positive elements. In \cite{SumitKaushik}, both  SQ  and SlerpSQ
have a similar effect in
interpolation, we choose SQ interpolation for fiber tracking.

\begin{figure}[htbp]
	\centering
	\subfloat[LogE]{\includegraphics[width=0.23\textwidth]{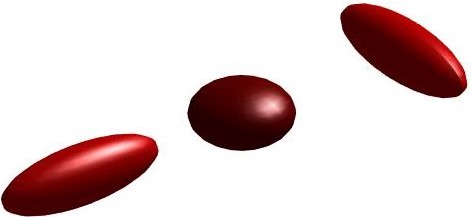}} \hspace{0.3cm}
	\subfloat[SQ]{\includegraphics[width=0.23\textwidth]{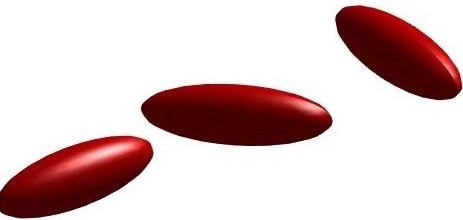}}\\
	\subfloat[Near the tensor at (12,13) location]{\includegraphics[width=0.155\textwidth]{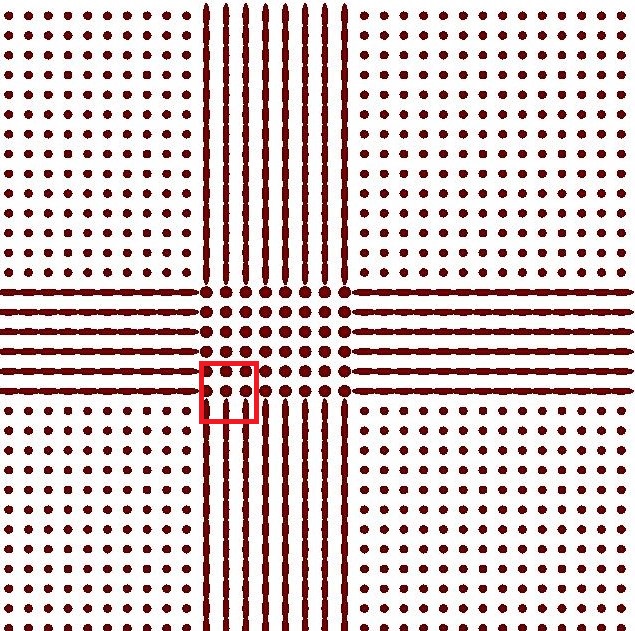}} \hspace{0.2cm}
	\subfloat[LogE interpolation]{\includegraphics[width=0.147\textwidth]{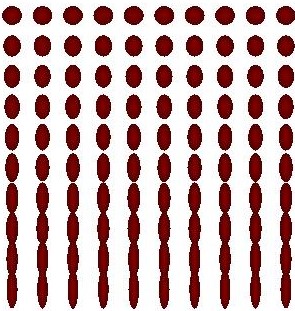}}\hspace{0.2cm}
	\subfloat[SQ interpolation]{\includegraphics[width=0.15\textwidth]{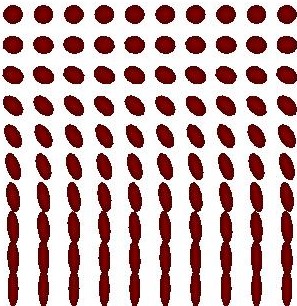}}
	\caption{Local Interpolation of 2nd order tensors in crossing fibers using LogE and SQ interpolation}
	\label{fig:interpolation using LogE}
\end{figure}

Fig \ref{fig:interpolation using LogE}(a) and (b) show the
spectral metric preserves anisotropy, which is crucial in fiber tracking
application.  In Fig \ref{fig:interpolation using LogE}(c), a tensor located
at (12,13), the rectangular section is considered.  The tensor at this
position is part of the vertical fiber and is underlying to the uniform
background and not cross-section with horizontal fiber.  Fig
\ref{fig:interpolation using LogE}(d) is interpolation in its neighborhood 
using the LogE metric and Fig \ref{fig:interpolation using LogE}(e) does the
same with respect to the SQ
metric.  The flow of interpolated tensors is more accurately captured in Fig
\ref{fig:interpolation using LogE}(e), which shows interpolation flows
towards the left.  Spectral metrics are known for robustness with respect to
noise, in segmentation of curved fibers, and presentation of anisotropy.

%\newpage
\section{Experiments and Results}
\subsection{Results on Synthetic Data}
For the experiments, we generate synthetic tensor fields with many configurations that have similar properties to many white matter tracts in the
brain. The synthetic images are simulated using a signal generated with b-value $1500 s/mm^{2}$ with a signal without gradient impulse $S_{0}=1$. Total of 81 gradient directions are chosen, which are uniformly distributed over the sphere. We use the adaptive kernel method to create fibers as detailed in \cite{adaptivbybarmutis}.

\begin{figure}[htbp]
\centering
\subfloat[Deep inverted U-shaped]{\includegraphics[width=0.22\textwidth]{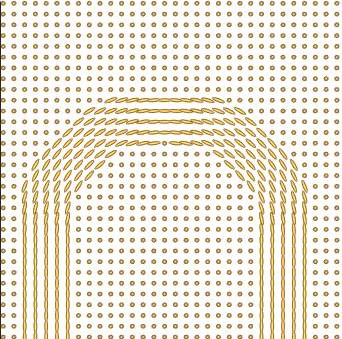}}\\
\subfloat[adjugate]{\includegraphics[width=0.23\textwidth]{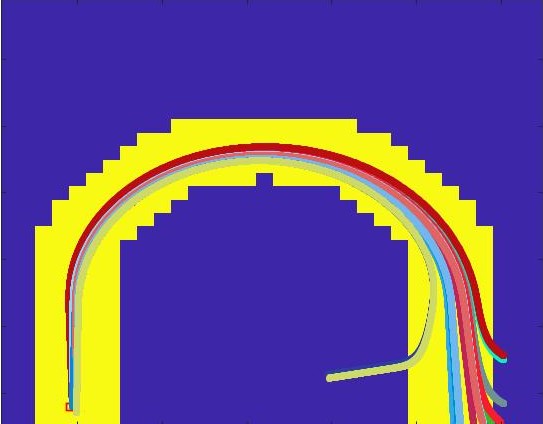}}\hspace{0.1cm}
\subfloat[$\beta$-scaled $p=2$]{\includegraphics[width=0.225\textwidth]{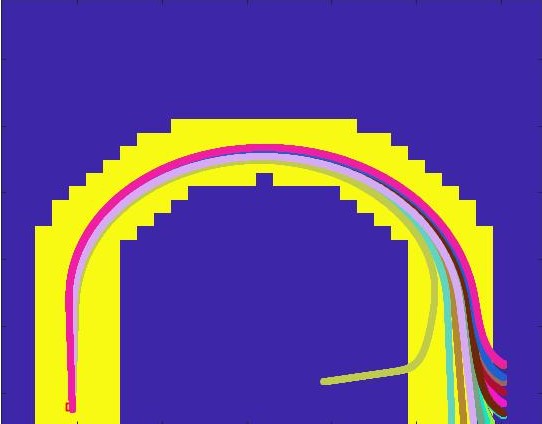}}\\
\subfloat[$\beta$-scaled $p=4$]{\includegraphics[width=0.225\textwidth]{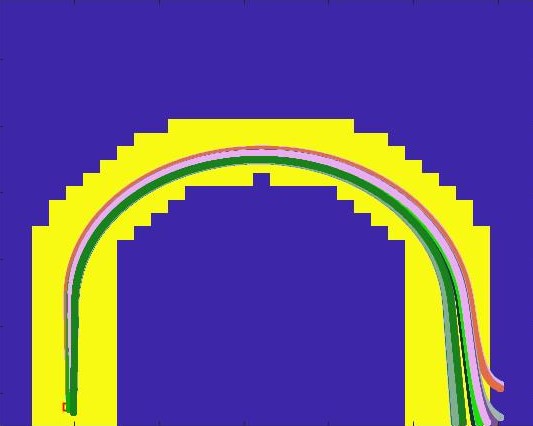}}\hspace{0.1cm}
\subfloat[$\beta$-scaled $p=6$]{\includegraphics[width=0.225\textwidth]{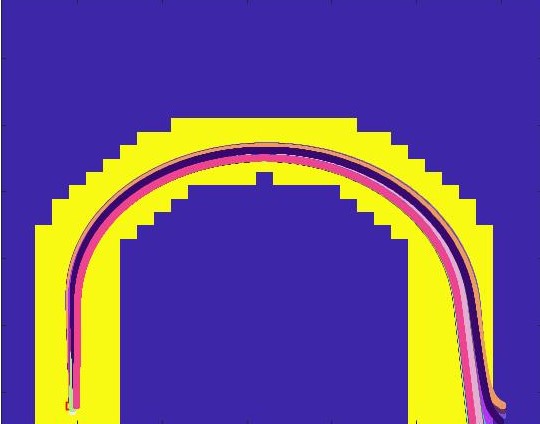}}
\caption{FA images corresponding to (a),comparison of ray-tracing with ODE solver using adjugate, and
$\beta$-scaled metric tensors from (b) - (e) respectively.}
\label{fig:comparision in inverted U}
\end{figure}

\begin{figure}[htbp]
	\centering
	\subfloat[4th order tensor]{\includegraphics[width=0.2\textwidth]{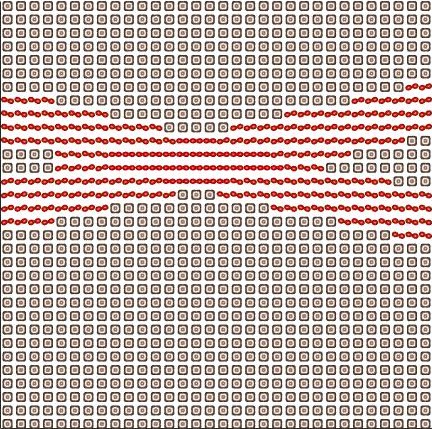}} \hspace{0cm}
	\subfloat[Diagonal sum]{\includegraphics[width=0.205\textwidth]{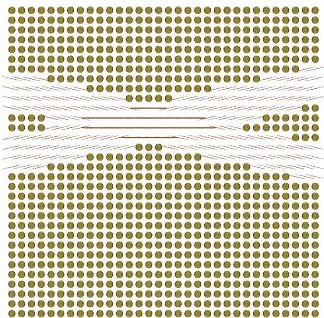}}\hspace{0cm}\\
	\subfloat[$\beta$-scaled tracking]{\includegraphics[width=0.2\textwidth]{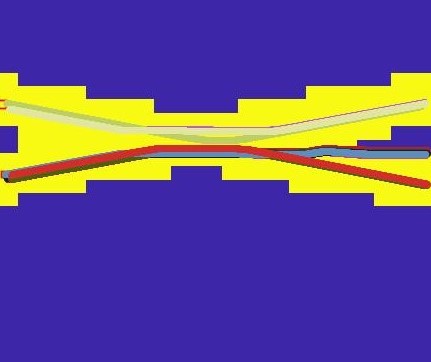}}\hspace{0cm}
	\centering
	\caption{Fiber tracing using $\beta$-scaled and DC  when two fibers are close.}
	\label{fig:sum of diagonal component}
\end{figure}

In Fig~\ref{fig:comparision in inverted U}, deep inverted U-shape is
considered with four points in the starting region and five shooting per
point.  It is visible that the adjugate metric fails to trace the fibers as
it approaches the target region.  $\beta$-scaled metric tensor geodesics
follow the fibers well, and higher power of $p$ produces smooth fibers and
increases fiber density.  Fig~\ref{fig:sum of diagonal component} shows
tractography result on the layer of the diagonal components where the two
fibers cross closely (cf.  Fig~\ref{fig:4th order Using DC}(a)).

\begin{figure}[htbp]
\centering
\subfloat[4th order inverted deep U shaped]{\includegraphics[width=0.22\textwidth]{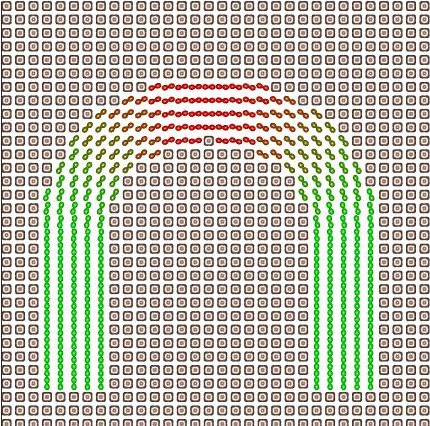}}\hspace{0.1cm}
\subfloat[4th order in Reflected S shape]{\includegraphics[width=0.22\textwidth]{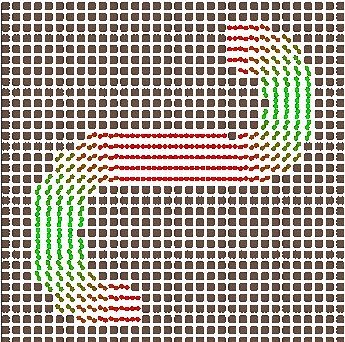}}\hspace{0.01cm}\\
\subfloat[Diagonal sum]{\includegraphics[width=0.22\textwidth]{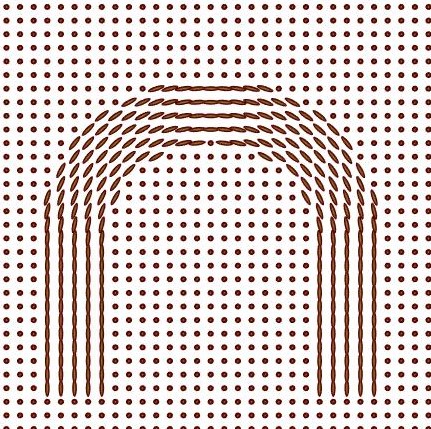}}\hspace{0.1cm}
\subfloat[Diagonal sum]{\includegraphics[width=0.22\textwidth]{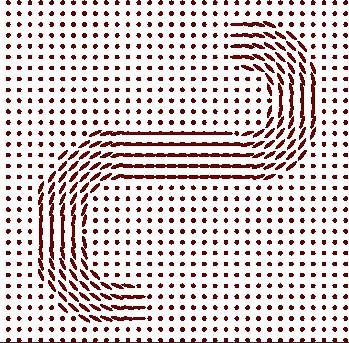}}\hspace{0.1cm}\\
	\subfloat[$\beta$-scaled tracing]{\includegraphics[width=0.235\textwidth]{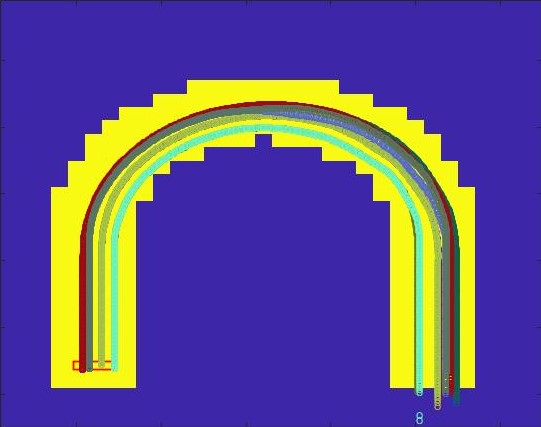}}\hspace{0.1cm}
	\subfloat[$\beta$-scaled tracing]{\includegraphics[width=0.235\textwidth]{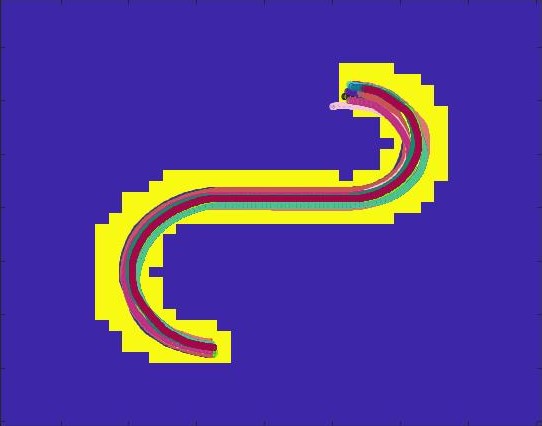}}
 \caption{Ray-tracing using Hybrid approach and diagonal sum 2nd order in high curvature fiber flows. Here we consider  with 10 points and 5 shoots per points inside the rectangular region.}.
	\label{fig:comparision using hybrid approach}
\end{figure}

%\newpage
Fig~\ref{fig:comparision using hybrid approach} shows that the hybrid approach can trace in high curvature fiber flows. On top of that, the 4th order tensor field image is shown with 2nd order tensor field obtained by sum of the diagonal components. This produces sharp images contrary to DTI.  \\
In Algorithm 1, the ODE solver method increases the
deviation of geodesic along the path. To overcome this
problem, we feedback the principal eigenvector direction
of the underlying interpolated tensor to ODE solver. This causes geodesic
deviation to disappear and leads to better performance under
the three metrics. We tested our method on different high
curvature fiber flows, as shown in Fig~\ref{fig:comparision using hybrid approach}.

\begin{figure}[htbp]
\centering
\subfloat[]{\includegraphics[width=0.22\textwidth]{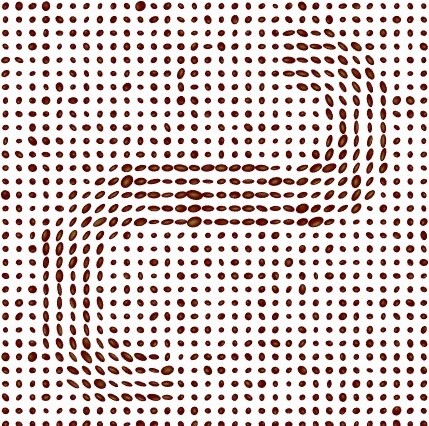}}\hspace{0.01cm}
\subfloat[]{\includegraphics[width=0.22\textwidth]{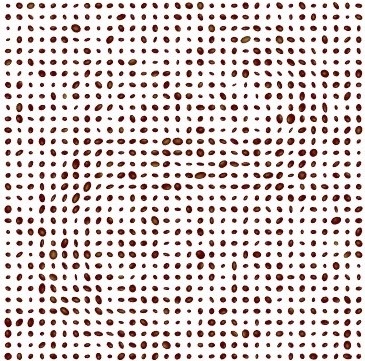}} \hspace{0.01cm}\\
\subfloat[]{\includegraphics[width=0.22\textwidth]{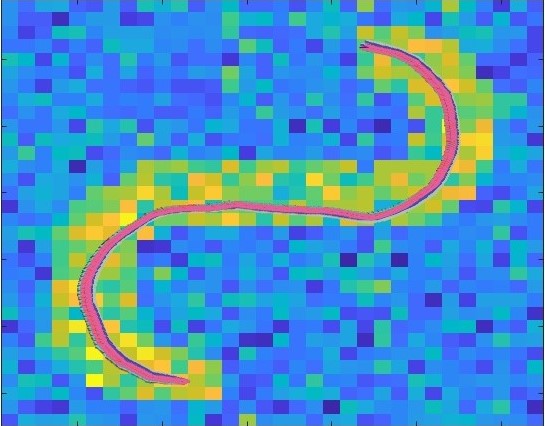}}\hspace{0.01cm}
\subfloat[]{\includegraphics[width=0.22\textwidth]{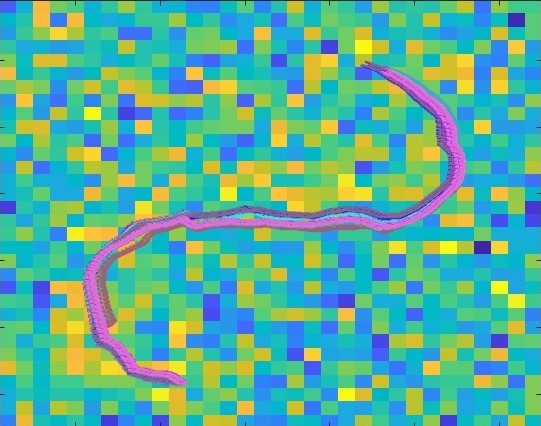}} \hspace{0.01cm}\\
\subfloat[]{\includegraphics[width=0.22\textwidth]{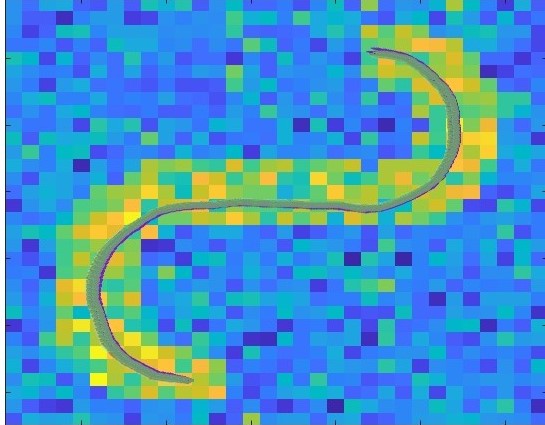}} \hspace{0.01cm}
\subfloat[]{\includegraphics[width=0.22\textwidth]{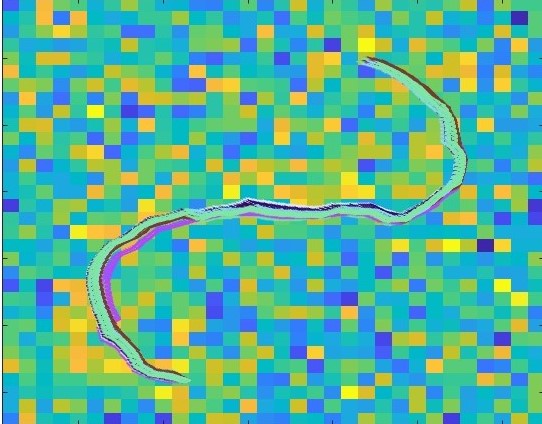}}
\caption{(a) Reflected S-shaped fiber with Riccian noise $0.25$, (b) Signal corrupted with Riccian noise $=0.30$, (c) Ray-tracing with Principal eigenvector direction using  adjugate and noise $0.25$ (d) Ray-tracing with Principal eigenvector direction using  adjugate and noise $0.30$, (e) Ray-tracing with principal eigenvector direction using
$\beta$-scaled metrics with $p=2$ and Riccian noise $0.25$, (f) Ray-tracing with principal eigenvector direction using
$\beta$-scaled metrics with $p=2$ and Riccian noise $0.30$}
\label{fig:comparision using noise}
\end{figure}

Fig~\ref{fig:comparision using noise} shows the S-shaped configuration corrupted with Riccian noise. The hybrid method acts robust and stable even in case  where the fiber is poorly visible Fig \ref{fig:comparision using noise}(b).
We show results of geodesic tracking on synthetic data for crossing fibers based on
$\beta$-scaled metric and diagonal component approach.

\begin{figure}[htbp]
\subfloat[4th order tensor ]{\includegraphics[width=0.23\textwidth]{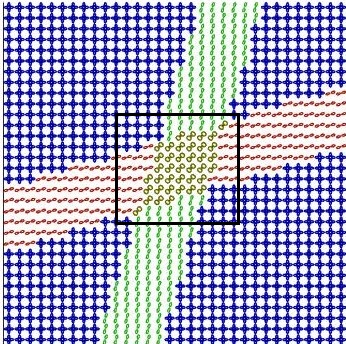}}\hspace{0cm}
\subfloat[Diagonal sum]{\includegraphics[width=0.24\textwidth]{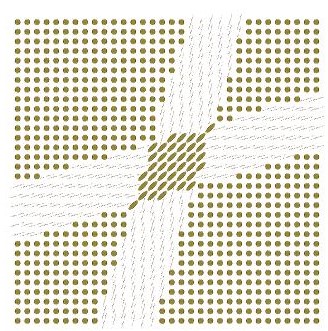}}\\ \hspace{0cm}
\subfloat[DC of rectangular section]{\includegraphics[width=0.21\textwidth]{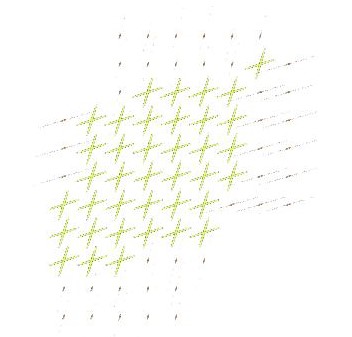}}\hspace{0.5cm}
 \subfloat[Fibers using $\beta$-scaled]{\includegraphics[width=0.24\textwidth]{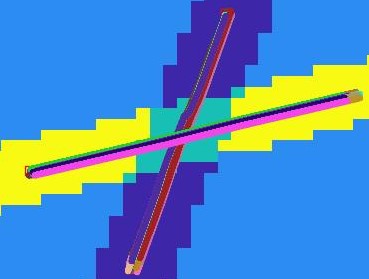}}	
 \caption{Fibers computed using the hybrid approach and diagonal sum at crossing area}
\label{fig:fiber tracking in crossing}
\end{figure}
\begin{figure}[htbp]
	\centering
	\subfloat[4th order]{\includegraphics[width=0.21\textwidth]{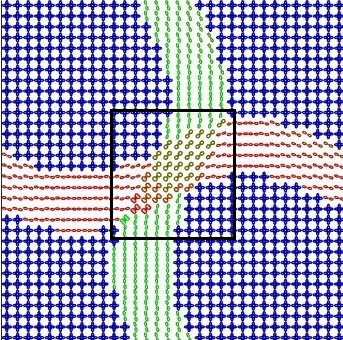}} \hspace{0.2cm}
	\subfloat[zoom section ]{\includegraphics[width=0.22\textwidth]{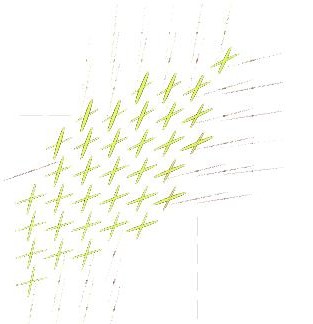}}\label{zoom _section}\hspace{0.1cm}
	\subfloat[first diagonal component]{\includegraphics[width=0.21\textwidth]{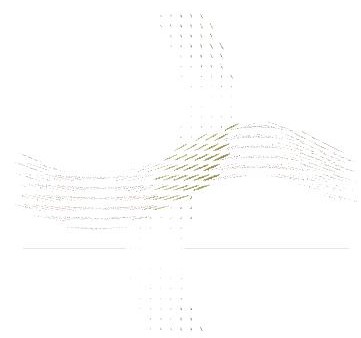}} \hspace{0.1cm}
	\subfloat[second diagonal component]{\includegraphics[width=0.21\textwidth]{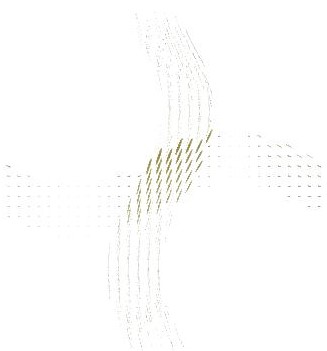}}\hspace{0.2cm}
	\subfloat[$\beta$-scaled  fiber tracking]{\includegraphics[width=0.24\textwidth]{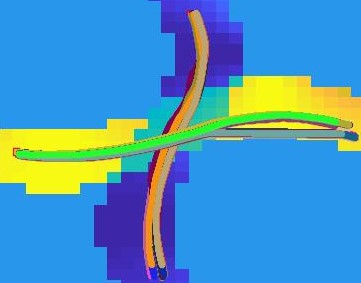}}
	\caption{This figure illustrate $\beta$-scaled fiber tracking can trace in high curvature fiber flow using the hybrid approach and diagonal sum at crossing area}
	\label{fig:fiber tracking sine curve}
\end{figure}
\begin{figure}[htbp]
	\centering
	\subfloat[Riccian noise 0.20]{\includegraphics[width=0.214\textwidth]{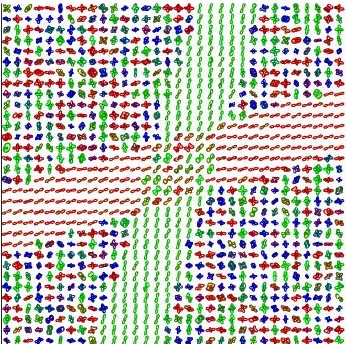}} \hspace{0cm}
	\subfloat[$\beta$-scaled tracking with noise 0.20]{\includegraphics[width=0.253\textwidth]{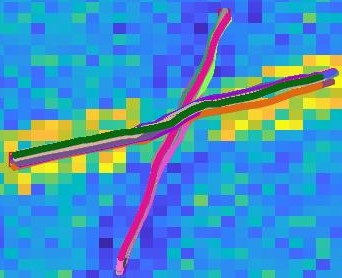}}\hspace{0.1cm}\\
\subfloat[4th order with noise 0.30]{\includegraphics[width=0.214\textwidth]{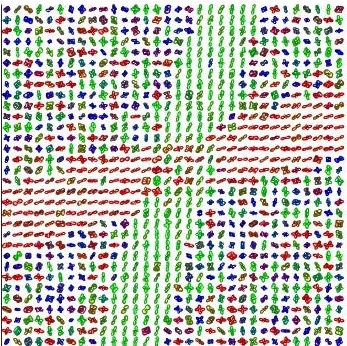}} \hspace{0cm}
	\subfloat[$\beta$-scaled tracking with noise 0.30]{\includegraphics[width=0.25\textwidth]{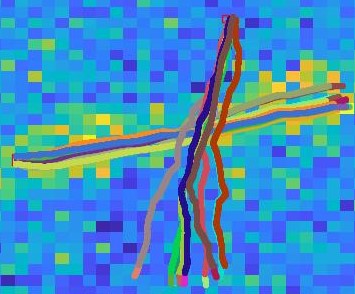}}\hspace{0cm}
	\caption{Fiber tracking in crossing fibers under noise}
	\label{fig:Fiber tracking under noise for crossing}
\end{figure}

Fig \ref{fig:fiber tracking in crossing} depicts two linear fibers
intersecting at small angle.  Fig \ref{fig:fiber tracking in crossing}(c)
shows the two components in the intersection region.  These components are
sharp and follows the trajectory of their individual single fibers.  Fig
\ref{fig:fiber tracking in crossing}(d) shows the $\beta$-scaled metric
tensor used for tracing the fiber bundle.

%\newpage
Fig \ref{fig:fiber tracking sine curve} has two curved fibers intersecting.  The Fig \ref{fig:fiber tracking sine curve}(b) shows the diagonal components in the intersection region. Fig \ref{fig:fiber tracking sine curve}(c) represent horizontally orientated regions of fiber curves  whereas \ref{fig:fiber tracking sine curve}(d) indicates the vertically oriented regions and \ref{fig:fiber tracking sine curve}(e) is the result of
$\beta$-scaled metric tensor tracking.  
In Fig \ref{fig:fiber tracking sine curve}(c and d), the diagonal components are able to align along with the correct running curved fibers.

%\newpage
In Fig~\ref{fig:Fiber tracking under noise for crossing}, the image shows two linear fibers crossing at sharp angles corrupted with noise. In this difficult case our method is able to reconstruct the fibers and behaves robust.

\subsection{Results on Real Data}

%\begin{flushleft}

Finally, we shortly comment on the results of our  tracking algorithm
applied to real images of the human brain.  

The DW-MRI image consists of
total size $114\times114\times70$ and each voxel is the size of
$2\times2\times2 \ mm^{3}$.  The real images are obtained by applying
gradient in 64 diffusion directions with diffusion weighting factor
$b=1500s/mm^{2}$ with single reference image (b = 0).  We have used
generalized logistic function (\ref{eq:A}) as activation to test on real
images.

\begin{figure}[htbp]
	\centering
	\subfloat[]{\includegraphics[width=0.21\textwidth]{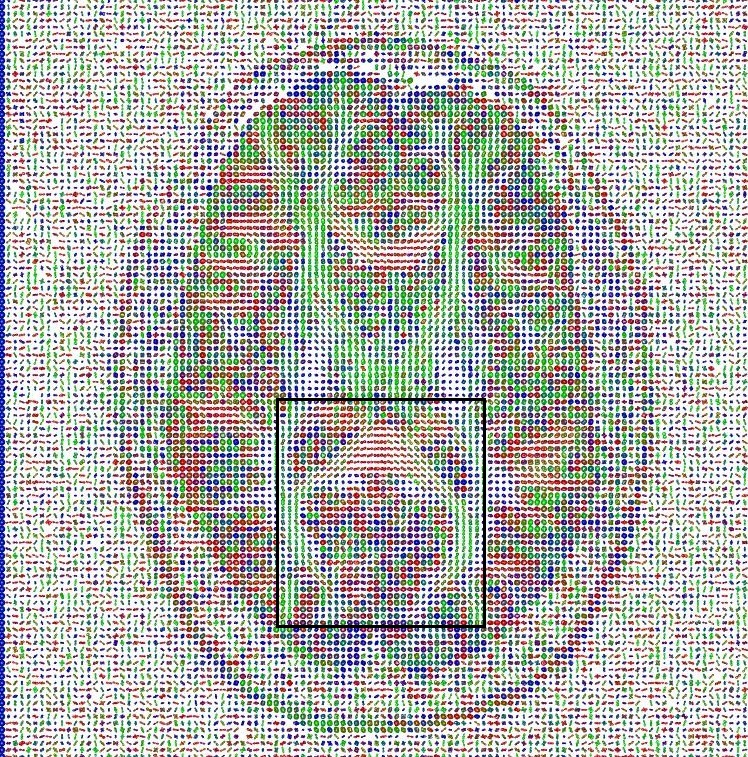}}\hspace{0cm}\\
	\subfloat[]{\includegraphics[width=0.21\textwidth]{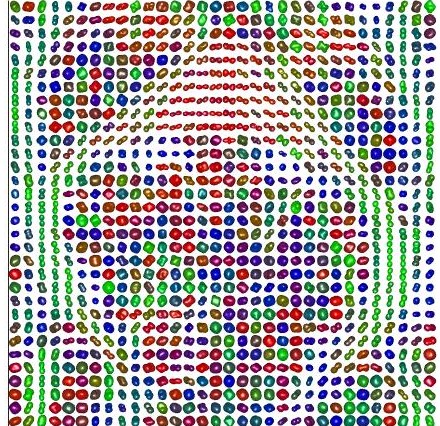}} \hspace{1mm}
	\subfloat[]{\includegraphics[width=0.256\textwidth]{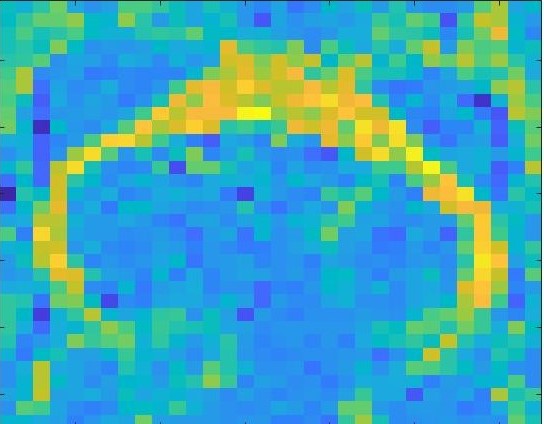}}	\hspace{0cm}
	\caption{(a) Real image (the selected rectangular section showing a section of Corpus Callosum, (b) 4th order tensor field in region of interest(black rectangular section in Fig (a)), (c) Fractional Anisotropy real Image}
	\label{fig:Real image with FA image}
\end{figure}

In Fig \ref{fig:Real image with FA image}, a rectangular section of the
Dorsal Longitudinal Fasciculus (DLF) is selected.  The Fig~\ref{fig:Real
image with FA image}(c) is corresponding to FA scalar image, while (b)
depicts the 4th order tensor.

%\end{flushleft}

\begin{figure}[htbp]
	\centering
	\subfloat[]{\includegraphics[width=0.24\textwidth]{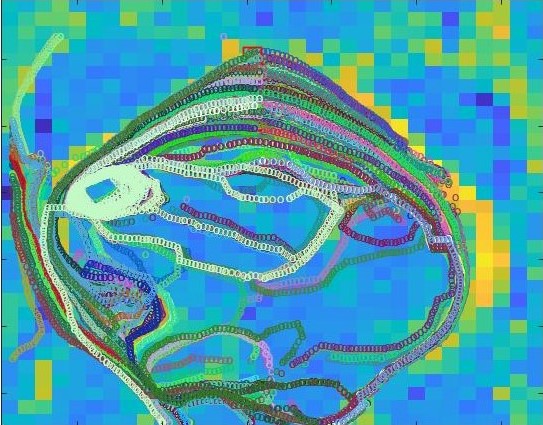}}\\
	\subfloat[]{\includegraphics[width=0.23\textwidth]{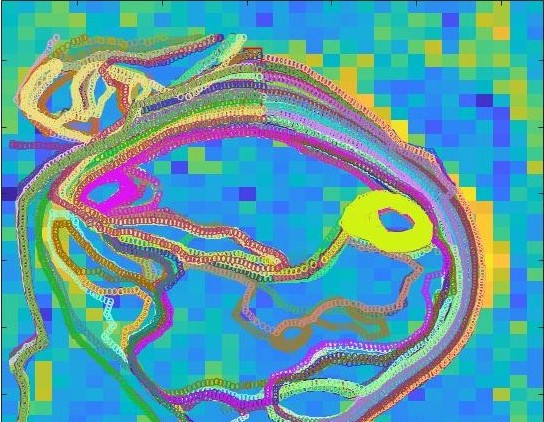}}\hspace{0.1cm}
	\subfloat[]{\includegraphics[width=0.23\textwidth]{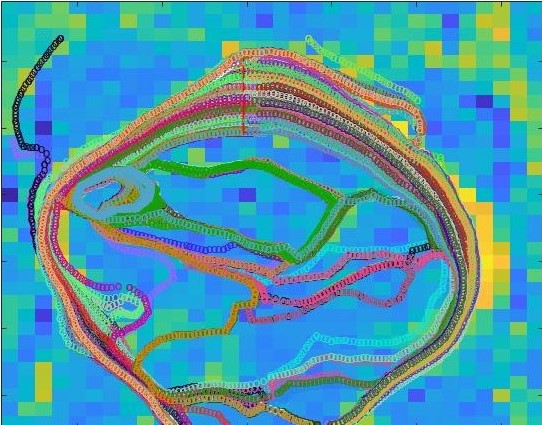}}	
	\caption{Fibers computed using the hybrid approach in Real Data with (a) $D^{-1}$, (b) Adjugate metric and (c)
	$\beta$-scaled metric with $p=2$}
	\label{fig:Tracking test on Real image}
\end{figure}

%\newpage
In Fig~\ref{fig:Tracking test on Real image}, 10 points were randomly picked
in rectangular section on the top of the structure.  Five geodesics are shot
per point in both directions.  The results indicate that most of the fibers
trace the white matter structure in all three cases.  Results under
$\beta$-scaled metric tensor produce smoother geodesics. 

The fiber crossing
resolution method proposed above was not used in this experiment.
%\newpage

\section{Discussion}
%\begin{flushleft}

We propose a new geodesic based tractography method by using a
$\beta$-scaled metric tensor.  This metric tensor is adapted according to
the inherent anisotropy property.  We have shown that the performance of
adapted metrics by means of sigmoid function as activation function composed
with the Hilbert anisotropy is better than performance of the classical
metric and it also performs better than adjugate metric for highly curved
fiber flows.  

To increase the accuracy of the tracking approach, we iterate
local geodesics tracing via Runge-Kutta ODE solver in the interpolated grid
of tensor data, initiated by the
principal eigenvector direction, called the hybrid
approach. 
%Here we also employ the spectral distance measures for the local
%interpolation.

Further, we propose to exploit the potential of using the diagonal
components of 4th order tensors, in particular for capturing crossing
fibers. 
%We resolve the fiber tracking in intersecting fibers region by
%using flattened 4th order tensors components.  
These diagonal components are second
order tensors lying in Riemannian symmetric space.  We have shown that they
have potential to effectively locate orientation distribution functions
(ODFs) maxima even at small angle intersections.

In future work, we will systematically use the spectral metric approach
for local interpolation.  The experiment discussed in section \ref{local
Interpolation} suggests to use spectral metric for local interpolation to
preserve anisotropy which is crucial for fiber reconstruction.

We also plan to employ the novel fiber crossing resolution approach in global
framework adapted for fiber tracking in real images.

%\end{flushleft}
%\newpage

\section*{Acknowledgments}
%\begin{flushleft}
The first three authors have been supported by the grant MUNI/A/0885/2019 of
Masaryk University, Jan Slov\'ak gratefully acknowledges support from the
Grant Agency of the Czech Republic, grant Nr.  GA20-11473S.  

We acknowledge
the core facility MAFIL supported by the Czech-BioImaging large RI project
(LM2018129 funded by MEYS CR) for their support with obtaining scientific
data presented in this paper.  

The second author acknowledges the support of
the OP VVV funded project "CZ.02.1.01/0.0/0.0/16\_019/0000765" Research
Center for Informatics, CTU in Prague.
%\end{flushleft}

%\bibliographystyle{unsrt}
%\renewcommand{bibfile}{References}
%\bibliographystyle{splncs}

  %\bibliographystyle{plain} % for sig-alternate
    %\bibliographystyle{unsrt}
    %\bibliographystyle{unsrtnat}
% \newpage
\bibliographystyle{IEEEtran}
\bibliography{bibfile}   

\end{document}